\documentclass[12pt]{article}
\usepackage[T1]{fontenc}
\usepackage[utf8]{inputenc}

\usepackage{moreverb}
\usepackage{a4wide}
\usepackage{authblk}
\usepackage[colorlinks,bookmarksopen,bookmarksnumbered,citecolor=red,urlcolor=red]{hyperref}
\usepackage{verbatim}
\usepackage{graphicx}
\usepackage{amsmath, amssymb}
\usepackage{multirow}
\usepackage{accents}
\usepackage{hhline}
\usepackage{amsmath,mathtools}
\usepackage{amsfonts}
\usepackage{amssymb}
\newtheorem{theorem}{Theorem}[section]
\newtheorem{corollary}[theorem]{Corollary}
\newtheorem{lemma}[theorem]{Lemma}
\newtheorem{remark}[theorem]{Remark}
\newtheorem{example}[theorem]{Example}
\newtheorem{definition}[theorem]{Definition}
\newtheorem{proposition}[theorem]{Proposition}

\newenvironment{proof}{\textbf{Proof:}}{\hfill$\Box$\\}

\newcommand{\R}{\mathbb{R}}

\newcommand{\RR}{\mathbb{R}}
\newcommand{\NN}{\mathbb{N}}

\newcommand{\calA}{{\cal A}}
\newcommand{\calL}{{\cal L}}
\newcommand{\bfu}{{\bf u}}
\newcommand{\bff}{{\bf f}}
\newcommand{\bPsi}{{\bf \Psi}}
\newcommand{\rba}{r^\beta_\alpha}

\title{Positive approximations of the inverse of fractional powers of 
SPD M-matrices}
\author[1]{S.~Harizanov\thanks{sharizanov@parallel.bas.bg}}
\author[1]{S.~Margenov\thanks{margenov@parallel.bas.bg}}
\affil[1]{{\small Institute of Information and Communication Technologies,\newline Bulgarian Academy of 
Sciences, Acad. G. Bonchev, bl. 25A, 1113 Sofia, Bulgaria}}

\begin{document}
\maketitle
\begin{abstract}
This study is motivated by the recent development in the fractional 
calculus and its applications. 
{
During last few years, several different techniques are  proposed to localize 
the nonlocal fractional diffusion operator. They are based on transformation 
of the original problem to a local elliptic or pseudoparabolic  problem, 
or to an integral representation of the solution,
thus increasing the dimension of the computational domain.
}
More recently, an alternative approach aimed at reducing the 
computational complexity was developed. The linear algebraic 
system $\cal A^\alpha \bf u=\bf f$, $0< \alpha <1$ is considered, 
where $\cal A$ is a 
{
properly normalized (scalded) symmetric and positive 
definite matrix obtained from finite element or finite difference
}
approximation of second order 
elliptic problems in 
{
$\Omega\subset\R^d$,
} 
$d=1,2,3$. The method is based on 
best uniform rational approximations (BURA) of the function  
$t^{\beta-\alpha}$ for $0 < t \le 1$ and natural $\beta$.

The maximum principles are among the major qualitative properties of
linear elliptic operators/PDEs. In many studies and applications,
it is important that such properties are preserved by the selected numerical
solution method. In this paper we present and analyze the properties
of positive approximations of $\cal A^{-\alpha}$ obtained by
{
the BURA technique. Sufficient conditions for positiveness are proven, 
complemented by sharp error estimates. The theoretical results are
supported by representative numerical tests.
}
\end{abstract}

\section{Introduction}\label{section1}
This work is inspired by the recent development in the fractional 
calculus and its various applications, i.e., to  Hamiltonian chaos, 
\cite{zaslavsky2002chaos}, 
anomalous diffusion in complex 
systems, \cite{bakunin2008turbulence}, long-range interaction in elastic 
deformations, \cite{silling2000reformulation}, nonlocal electromagnetic fluid 
flows, \cite{mccay1981theory}, image processing, \cite{gilboa2008nonlocal}.
A more recent impressive examples of anomalous diffusion models in chemical 
engineering are provided in \cite{metzler2014anomalous}.
Such kind of applications lead to fractional order partial differential 
equations that involve in general non-symmetric elliptic operators see, e.g.
\cite{KilbasSrivastavaTrujillo:2006}. An important subclass of this topic are 
the fractional powers of self-adjoint elliptic operators, which are nonlocal 
but 
{
self-adjoint. 
}
In particular, the fractional Laplacian \cite{Pozrikidis16}
describes an unusual diffusion process associated with random excursions.
In general, the parabolic equations with fractional derivatives in time are 
associated with sub-diffusion, while the fractional elliptic operators are
related to super-diffusion. 

Let us consider the elliptic boundary value problem in a weak form:  
find  $u \in V$ such that
\begin{equation}\label{eqn:weak}
a(u,v) := \int_\Omega \left ({\bf a}({x}) \nabla u({x}) \cdot  \nabla v({x}) + q(x) 
\right )d{x} =\int_\Omega f({x}) v({x}) d {x},~~~ \forall v \in V,
\end{equation}
where
$$V:= \{ v \in H^1(\Omega): ~~v({x})=0 ~~\mbox{on} ~~\Gamma_D \},
$$
$\Gamma = \partial \Omega$, and
$\Gamma = \bar \Gamma_D \cup \bar \Gamma_N$.
We assume that $\Gamma_D$ has positive measure, $q({x}) \ge 0$ in 
$\Omega$, and ${\bf a}({x})$ 
is an SPD $d \times d$ matrix, uniformly bounded in $\Omega$, i.e.,
\begin{equation}\label{ax}
c \Vert {\bf z} \Vert^2  \leq
{\bf z}^T {\bf a}({x}) \, {\bf z} \leq
C \Vert{\bf z} \Vert^2 \quad \forall {\bf z} \in\R^{d}, \forall {x} \in \Omega,
\end{equation}
for some positive constants $c$ and $C$. Also, $\Omega$
is a polygonal domain in 
$\R^d$, $d\in\{1,2,3\}$, and $f({x})$ is a
given Lebesgue integrable function on $\Omega$ that belongs to the space $L_2(\Omega)$.
Further, the case when ${\bf a}({x})$ does not depend on $x$ is referred 
to as problem in  homogeneous media,
while the general case models processes in non homogeneous media.
The bilinear form $a(\cdot, \cdot)$ defines a linear operator 
${\cal L}: V \to V^*$ with $V^*$ being the dual of $V$. Namely, 
for all $u,v \in V$ $a(u,v) := \langle {\cal L} u,  v \rangle$,
where $\langle \cdot, \cdot \rangle$ is the pairing between $V$ and $V^*$.

One possible way to introduce  ${\cal L}^{\alpha}$, $0<\alpha<1$, is through its 
spectral decomposition, i.e.
\begin{equation}\label{SpDec}
{\cal L}^{\alpha} u(x) = \sum_{i=1}^\infty \lambda_i^{\alpha} c_i \psi_i(x), \quad
\mbox{where} \quad
u(x) = \sum_{i=1}^\infty  c_i \psi_i(x).
\end{equation}
Here
$\{ \psi_i(x) \}_{i=1}^\infty$ are the eigenfunctions of $\cal L$,
orthonormal in $L_2$-inner product and  $\{ \lambda_i \}_{i=1}^\infty$
are the corresponding positive real eigenvalues.
This definition generalizes the concept of equally weighted left and right 
Riemann-Liouville fractional derivative, defined in one space dimension, to the 
multidimensional case. There is still ongoing research about the relations of 
the different definitions and their applications, see, e.g. \cite{bates2006nonlocal}.

The numerical solution of nonlocal problems is rather expensive. 
The following three approaches (A1 - A3) are based on transformation 
of the original problem 
\begin{equation}\label{NonLP}
{\cal L}^{\alpha} u = f
\end{equation}
to a local elliptic or pseudo-parabolic problem, or on integral 
representation of the solution, thus increasing the dimension of 
the original computational domain. 

The Poisson problem is considered in the related papers refereed bellow, i.e. 
$$
a(u,v) := \int_\Omega \nabla u({x}) \cdot  \nabla v({x})  d{x}.
$$

\begin{itemize}
\item[{\bf A1}] 
{
Extension to a mixed boundary value problem in the 
semi-infinite cylinder $C=\Omega \times \mathbb R_+ \subset \mathbb R^{d+1}$

A “Neumann to Dirichlet” map is used in \cite{CNOSaldago_2016}. 
Then, the solution of fractional Laplacian problem is obtained by 
$u(x) = v(x, 0)$ 
where $v : \Omega\times\mathbb R_+ \rightarrow \mathbb R$ is a solution of the equation
$$
-div\left ( y^{1-2\alpha}\nabla v(x,y)\right ) = 0, ~~~ 
(x,y)\in \Omega\times \mathbb R_+ ,
$$
where $v(\cdot,y)$ satisfies the boundary conditions 
of (\ref {eqn:weak}) $\forall y\in \mathbb R_+$,
$$
\lim_{y\rightarrow\infty} v(x,y) = 0, ~~~ x\in\Omega ,
$$
as well as
$$
\lim_{y\rightarrow 0^+} \left (- y^{1-2\alpha} v_y(x,y)\right ) = f(x), 
~~~ x\in\Omega.
$$
It is shown that the variational formulation of this equation is 
well posed in the related weighted Sobolev space. The finite element 
approximation uses the rapid decay of the solution $v(x,y)$ in the 
$y$ direction, thus enabling truncation of the semi-infinite
cylinder to a bounded domain of modest size. The proposed multilevel 
method is based on the Xu-Zikatanov identity \cite{XuZ2002}. 
The numerical tests for $\Omega=(0,1)$ and $\Omega=(0,1)^2$  
confirm the theoretical estimates of almost optimal computational complexity. 
}

\item[\bf A2]
{
Transformation to a pseudo-parabolic problem  

The problem (\ref{eqn:weak}) is considered in \cite{Vabishchevich14,Vabishchevich15} assuming the boundary 
condition 
$$
a(x){\frac{\partial u}{\partial n}} + \mu (x) u = 0, ~~~ x\in \partial\Omega,
$$
which ensures ${\cal L} = 
{\cal L}^* \ge \delta {\cal I}$, $\delta >0$. Then the solution of 
fractional power diffusion problem $u$ can be 
{found
}
 as
$$
u(x)=w(x,1), ~~~ w(x,0)=\delta^{-\alpha}f,
$$
where $w(x,t), 0<t<1$, is the solution of pseudo-parabolic equation 
$$
(t{\cal D} + \delta {\cal I} ) {\frac{dw}{dt}} + \alpha {\cal D}w = 0,
$$
and ${\cal D} = {\cal L} - \delta {\cal I} \ge 0$. Stability conditions 
are obtained for the fully discrete schemes under consideration.
A further development of this approach is presented in \cite{LV17} where
the case of fractional order boundary conditions is studied.
}

\item[\bf A3]
{
Integral representation of the solution
 
The following representation of the solution of 
(\ref {eqn:weak}) 
is used
in \cite{BP15}:
$$
{\cal L}^{-\alpha} = \frac{2\sin (\pi\alpha)}{\pi} \int_0^\infty  t^{2\alpha - 1} 
\left ({\cal I} +t^2 {\cal L}\right )^{-1} dt 
$$
Among others, the authors introduce an exponentially convergent  
quadrature scheme. Then, the approximate solution of $u$ only involves evaluations 
of $({\cal I}+t_i {\cal A})^{-1} f$, where $t_i\in (0,\infty)$ is 
related to the current quadrature node, and where $\cal I$ and ${\cal A}$ stand 
for the identity and the finite element stiffness matrix corresponding
to the Laplacian. The computational complexity of the method depends 
on the number of quadrature nodes. For instance, the presented analysis 
shows that approximately 50 auxiliary linear systems have to be solved to get 
accuracy of the quadrature scheme of order $O(10^{-5})$ for 
$\alpha \in \{0.25, 0.5, 0.75\}$.  
{
A further development of this approach is available in \cite{BP16}, where the theoretical analysis is
extended to  the class of regularly accretive operators.
}
}
\end{itemize}

An alternative approach is applied in \cite{HLMMV16} where a class of optimal 
solvers for linear systems with fractional power of symmetric and positive 
definite (SPD) matrices is proposed.
Let $\calA \in \R^{N \times N}$ be a 
{normalized} SPD matrix generated
by a finite element 
{or finite difference }
approximation of some self-adjoint elliptic problem. 
An efficient method for 
solving algebraic systems of linear equations involving fractional powers
of the matrix $\calA$ is considered, namely for solving the system
\begin{equation}\label{eq:fal}
\calA^\alpha {\bfu} ={\bff},  \quad \mbox{ where} \quad  0 < \alpha <1.
\end{equation}
The fractional power of SPD matrix $\calA$, similarly to the infinite 
dimensional counterpart $\calL$, is expressed through the spectral representation
of $\bfu$ through the eigenvalues and eigenvectors $\{ (\Lambda_i, \bPsi_i) \}_{i=1}^N $
of $\calA$, assuming that the eigenvectors are $l_2$-orthonormal, i.e.
$\bPsi_i^T \bPsi_j = \delta_{ij}$  and $\Lambda_1 \le \Lambda_2 \le \dots \Lambda_N
{\le1}$.
Then $\calA={\cal W}{\cal D}{\cal W}^T$, $\calA^\alpha = 
{\cal W}{\cal D}^\alpha {\cal W}^T$,
where the $N \times N$ matrices  ${\cal W}$ and ${\cal D}$ are
defined as ${\cal W}=(\bPsi_1^T, \bPsi_2^T, ..., \bPsi_N^T)$ and
${\cal D}=diag(\Lambda_1, \dots, \Lambda_N)$, 
$ \calA^{-\alpha}= {\cal W} {\cal D}^{-\alpha} {\cal W}^T$, and the solution of
$\calA^\alpha \bfu =\bff$ can be expressed as
\begin{equation}\label{eqn:exacalg}
{\bfu} = \calA^{-\alpha} \bff = {\cal W} {\cal D}^{-\alpha} {\cal W}^T {\bf f}.
\end{equation}
Instead of the system (\ref{eq:fal}), one can solve the equivalent system
$\calA^{\alpha  - \beta} \bfu= \calA^{-\beta} \bff
{:=\mathbf F} $ with $ \beta \ge 1$ an integer. 
Then  the idea is to approximately evaluate $\calA^{\beta -\alpha} 
{\mathbf F} $ using a 
set of equations involving inversion of $\calA$ and $\calA - d_j I $, 
for $j=
{1}, \dots, k$. The integer parameter 
$k \ge 1$ is the number of partial fractions of the best uniform rational approximation 
(BURA) $r^\beta_\alpha(t)$  of $t^{\beta -\alpha}$ on the interval $(0,1]$. One can observe
that the algorithm of \cite{BP15}, see A3, can be viewed as a particular rational
approximation of $\calA^{-\alpha}$. It is also important, 
{
that in certain sense the results from 
\cite{HLMMV16} are more general, and are applicable to a wider class of sparse SPD matrices.
}

Assuming that $\calA$ is a large-scale matrix, the computational complexity
of the discussed methods for numerical solution of fractional diffusion problems
is substantially high. Then, the parallel implementation of such methods for real
life problems is an unavoidable topic. 
{
In this context, there are some serious advantages of the last two approaches, see e.g., in \cite{CSM16}.
}

The maximum principles are among the major qualitative properties of the 
elliptic or parabolic operators/PDEs. In general, to solve 
PDEs we use some numerical method, and it is a natural requirement that 
such qualitative properties are preserved on the discrete level. Most of 
the studies which deal with such topics give sufficient conditions for the 
discretization parameters in order to guarantee the certain maximum principle. 
It is easily see, that under certain such assumptions, 
the solutions of 
{
(A1) and (A3) 
} 
satisfy certain maximum principle. In this paper, we study
positive approximations of $\calA^{-\alpha}$, obtained by BURA technique
introduced in \cite{HLMMV16}, under rather general assumptions for 
{
the normalized SPD matrix $\calA$.
}   

The rest of the paper is organized as follows. In Section \ref{section2} 
we provide a brief introduction to the topic of monotone matrices including 
some basic properties of the M-matrices and their relations to FEM 
discretization of elliptic PDEs. Sufficient conditions for positive 
approximations of the inverse of a given SPD matrix,
based on BURA technique are presented in Section \ref{section3}. The analysis in
Section \ref{section4} is devoted to a class of best rational approximations of 
$\calA^{-\alpha}$, that satisfy such sufficient conditions. 
{Sharp error estimates for a class of BURA approximations are also included 
in this section. Some numerical tests and short concluding remarks are given at 
the end.}

\section{Monotone matrices and SPD M-matrices}\label{section2} 

The maximum principles are some of the most useful properties used to solve
a wide range of problems in the PDEs. For instance, their use is often essential 
to study the uniqueness and necessary conditions of solvability, approximation 
{
and boundedness
} 
 of the solution, as well as, for quantities of physical 
interest like maximum stress, torsional stiffness, electrostatic capacity, 
charge density etc. 
Under certain regularity conditions, a classical maximum principle for elliptic 
problems reads as follows. Suppose that $\calL u\ge0$ in $\Omega$, 
then a nonnegative maximum is attained at the boundary $\partial\Omega$.
Let us assume additionally that $u\ge 0$ in $\partial\Omega$. Then the positivity
preserving property holds, that is, $u(x)>0$ for $x\in\Omega$ or $u\equiv 0$. 

To solve PDEs we use some numerical methods, and it is a natural requirement 
that such qualitative properties are preserved on the discrete level. Most of 
the papers which deal with this topic give sufficient conditions for the 
discretization parameters in order to guarantee the certain maximum principle.
For instance, when FEM is applied, the related results are usually described in 
terms of properties of the related mass and stiffness matrices. 
{\definition {A real square matrix $\calA$ is called monotone if for all real vectors
$v$, $\calA {\bf v}\ge 0$ implies ${\bf v}\ge 0$, where $\ge$ is in element-wise sense.}}

\vspace{10pt}
\noindent
The next property is sometimes used as an alternative definition.
{\proposition {Let $\calA$ be a real square matrix. $\calA$ is monotone if and only if 
$\calA^{-1}\ge 0$.}}

{\definition 
{The class of Z-matrices are those matrices whose off-diagonal entries are less 
than or equal to zero. Let $\calA$ be a 
{$N\times N$} real Z-matrix, then $\calA$ 
is a non-singular M-matrix if every real eigenvalue of $\calA$ is positive. 
A symmetric M-matrix is sometimes called a Stieltjes matrix.}}

\vspace{10pt}
\noindent
The M-matrices are among most often used monotone matrices. They arise 
naturally in some discretizations of elliptic operators.  

Let $\calA$ be an SPD matrix obtained after FEM approximation of
(\ref{eqn:weak}) by linear triangle elements. Let us assume also 
that $\Omega\subset\R^2$ is discretized by a nonobtuse triangle 
mesh $\tau_h$, and the coefficients $a(x)=a_e$ and $q(x)=q_e$ are 
piecewise constants on the triangles $e\in\tau_h$. Then $\calA$ can be
assembled by the element matrices 
$$
\calA_e = {\cal K}_e + {\cal M}_e,
$$  
where ${\cal K}_e$ is the element stiffness matrix and 
${\cal M}_e$ is the element mass matrix.  
Subject to a scaling factor, the off diagonal elements of the 
symmetric and positive semidefinite matrix ${\cal K}_e$ are equal 
to some of $-\cot (\theta_i)\le 0$, $i=1,2,3$, where 
$0<\theta_i\le\pi/2$ are the nonobtuse angles of the triangle $e$. 
Imposing the boundary conditions we get that the global stiffness 
matrix is SPD M-matrix. The element mass matrix is positive 
diagonal matrix if a proper quadrature formula is applied.
A similar result can be obtained by the standard diagonalization 
known as {\it lumping the mass}. Then, the global mass matrix is 
positive diagonal matrix and $\calA$ is SPD M-matrix. A more
general considerations of this kind are available in \cite{KraMar2009} 
including the case of coefficient anisotropy as well as the 
nonconforming linear finite elements. Similar representations 
of the element mass and stiffness matrices are derived
if $\Omega\subset\R^3$ is discretized by a nonobtuse tetrahedral 
mesh $\tau_h$. We get again that $\calA$ is again SPD M-matrix, 
see e.g., \cite{KosMar2009}.

{\remark {
Not all monotone matrices are M-matrices, and the sum of two 
monotone matrices is not always monotone. The next examples 
prove these statements. 
}} 

{\example
~~\\
\begin{itemize}
\item[{\it E1.1}] 
{
$\calA_1 = \left ( 
\begin{matrix} 
-1 & 3\\
2 & -4
\end{matrix}
\right)$ 
is not M-matrix, but
$\calA_1^{-1} = {\frac{1}{2}} \left (\begin{matrix}
4 & 3\\
2 & 1
\end{matrix} \right )
$ and therefore $\calA_1$ is monotone.}
\item[{\it E1.2}]
{
$\calA_2=\calA_1+6{\cal I}$ is a sum of two monotone matrices, but
$\calA_2^{-1} =  {\frac{1}{4}} \left ( 
\begin{matrix}
2 & -3\\
-2 & 0
\end{matrix} \right ),
$ and therefore $\calA_2$ is not monotone.
}
\end{itemize}
}

\vspace{6pt}
In what follows, we study positive approximations of 
$\calA^{-\alpha}$ for 
{a given normalized} 
SPD M-matrix $\calA$. It follows straightforwardly that the 
inverse of each such approximation will approximate
$\calA^\alpha$ in the class of  monotone functions.

\section{Positive approximations of the inverse of SPD 
M-matrices}\label{section3}
Explicit computation and memory storage of the fractional power $\calA^{\alpha}$ in \eqref{eq:fal} for large-scale 
problems is expensive and impractical. Even when $\calA$ is sparse, $\calA^\alpha$ is typically dense. Therefore, we 
study possible positive approximations of the action of $\calA^{-\alpha}$ based solely on the information of $\calA$. 
For this purpose, we consider the class of rational functions
$$
\mathcal R(m,k):=\{P_m/Q_k\,:\,P_m\in{\mathcal P}_m,Q_k \in {\mathcal P}_k\},
$$ 
fix a positive integer $\beta$, and search for an appropriate candidate $r$ in it, that approximates well 
the univariate function $t^{\beta-\alpha}$ on the unit interval $[0,1]$. Note that, due to the normalization of 
$\calA$, this interval covers the spectrum of $\calA^\alpha$, 

{\definition {Let $\alpha\in(0,1)$, and $\beta,m,k\in\NN\setminus\{0\}$. The minimizer $r^\beta_\alpha\in\mathcal 
R(m,k)$ of the problem
\begin{equation}\label{eq:BURA}
\min_{ r \in \mathcal R(m,k)}  \max_{t \in [0,1]} \left | t^{\beta-\alpha} - r(t)  \right |,
\end{equation}
will be called $\beta$-Best Uniform Rational Approximation ($\beta$-BURA). Its error will be denoted by
$$
E_\alpha(m,k;\beta):=\max_{t \in [0,1]} \left |t^{\beta-\alpha} - r^\beta_\alpha(t)\right |.
$$
}}
Based on classical Spectral Theory arguments (see \cite[Theorem 2.1]{HLMMV16}), the univariate approximation error 
$E_\alpha(m,k;\beta)$ is an upper bound for the multivariate relative error 
$\|\calA^{-\beta}\rba(\calA)\bff-\calA^{-\alpha}\bff\|_{\calA^{\gamma+\beta}}/\|\bff\|_{\calA^{\gamma-\beta}}$ for the 
corresponding matrix-valued BURA approximation of the exact solution $\bfu$ in \eqref{eqn:exacalg}. Here, 
$\gamma\in\RR$ can be arbitrary, 
{
and the Krylov norms are defined via standard energy dot product, i.e. 
$\|\bff\|^2_{\calA^{\gamma-\beta}}=\langle 
\calA^{\gamma-\beta}\bff,\bff\rangle$.
}

{\proposition\label{prop2} {Let $\calA\in\RR^{N\times N}$ be an SPD matrix with eigenvalues 
$0<\Lambda_1\le\Lambda_2\le\dots\le\Lambda_N\le1$. Let $\rba$ be the $\beta$-BURA for given $\alpha,\beta,m,k$. Then, 
\begin{equation}\label{eq:BURA error}
\|\calA^{-\beta}\rba(\calA)\bff-\calA^{-\alpha}\bff\|_{\calA^{\gamma+\beta}}\le 
E_\alpha(m,k;\beta)\|\bff\|_{\calA^{\gamma-\beta}},\qquad\forall \gamma\in\RR,\; \forall\bff\in\RR^N.
\end{equation}
}}
For the practical computation of $\calA^{-\beta}\rba(\calA)\bff$ we use the partial fraction decomposition of 
$t^{-\beta}\rba(t)$, which is of the form
\begin{equation}\label{eq:FactFull}
t^{-\beta}r^\beta_\alpha(t) =
\sum_{j=0}^{m-k-\beta} b_j\, t^{j}+
\sum_{j=1}^{\beta}\frac{c_{0,j}}{t^j}    +\sum_{j=1}^{k} \frac{c_j}{t - d_j}, 
\end{equation}
provided all $\rba$ has no complex poles and the real ones $\{d_j\}_1^k$ are all of multiplicity 1. Later, we will see 
that for $\beta=1$ and $m=k$ the above assumption on the poles of $\rba$ holds true for any $\alpha\in(0,1)$. 
Furthermore, in all our numerical experiments with various $\beta,m,k$ and $\alpha\in\{0.25,0.5,0.75\}$ the assumption 
always remains valid. Hence, it does not seem to restrict the application range of the proposed method. On the other 
hand, under \eqref{eq:FactFull} the approximate solution 
\begin{equation}\label{eq:Rgeneral}
\bfu_r:=\calA^{-\beta}r^\beta_\alpha(\calA) \bff   =
\sum_{j=0}^{m-k-\beta}b_i\calA^j\bff +
\sum_{j=1}^{\beta}c_{0,j}\calA^{-j}\bff  +
\sum_{j=1}^{k} {c_j}{(\calA - d_j I)^{-1}} \bff 
\end{equation}
of $\bfu$ can be efficiently numerically computed via solving several linear systems, that involve $\calA$ and its 
diagonal variations $\calA-d_j I$, for $j=1,\dots,k$.

{\definition {A real symmetric matrix $\calA^{-1}$ is said to be doubly nonnegative if it is both positive 
definite, and entrywise nonnegative.}}

Our 
{first} goal is to analyze under what conditions on the coefficients and the poles in \eqref{eq:FactFull}, the matrix 
$\calA^{-\beta}r^\beta_\alpha(\calA)$ remains doubly nonnegative. Clearly the matrix is symmetric whenever $\calA$ is, 
so the main investigations are on assuring $\calA^{-\beta}r^\beta_\alpha(\calA)\ge0$. The following proposition 
contains sufficient conditions for positivity.

{\proposition\label{prop3} {If $\calA$ is a normalized SPD M-matrix, $m<k+\beta$, $c_0\ge0$, $c\ge 0$, 
{and ${\mathbf d}<0$ 
(entrywise)}, then 
$\calA^{-\beta}r^\beta_\alpha(\calA)$ in \eqref{eq:Rgeneral} is doubly nonnegative.}}
\begin{proof}
Since $m<k+\beta$, equation \eqref{eq:Rgeneral} is simplified to  
$$\calA^{-\beta}r^\beta_\alpha(\calA)=\sum_{j=1}^{\beta}c_{0,j}\calA^{-j} + \sum_{j=1}^{k} {c_j}{(\calA - d_j 
I)^{-1}}.$$
For every $j=1,\dots,k$, the matrix $\calA-d_j I$ is an SPD M-matrix, as $d_j<0$ and the 
diagonal elements increase their values,  
{i.e. become stronger dominant.} 
Hence, $(\calA-d_j I)^{-1}\ge0$.
 We have $\calA^{-1}\ge0$, thus $\calA^{-j}=(\calA^{-1})^j\ge0$, $j=1,\dots,\beta$, as each entry of $\calA^{-j}$ is a 
sum of nonnegative summands. Finally, a linear combination of positively scaled doubly nonnegative matrices is also a 
doubly nonnegative matrix. 
\end{proof}

Note that, when applying pure polynomial approximation techniques for $t^{-\alpha}$ on $[\Lambda_1,1]$ like in 
\cite{HMMV2016}, there is practically no chance to come up with a positive approximation of $\calA^{-\alpha}$. First of 
all, such an approximant is a linear combination of positive degrees of $\calA$ and in particular $\calA$ itself 
appears with a nonzero coefficient. This matrix has non-positive off-diagonal entries. Furthermore, it was numerically 
observed that the coefficient sequence in the linear combination is sign alternating. Hence, the proposed $\beta$-BURA 
approach seems the right and most natural tool for constructing positive approximations of $\calA^{-\alpha}$, or 
alternatively, monotone approximations of $\calA^\alpha$. Another disadvantage of the former approach is the 
restriction on $\Lambda_1$ to be well-separated from zero, which is also a restriction on the condition number of 
$\calA$.

\section{Analysis of a class of best rational approximations of fractional power 
of SPD M-matrices}\label{section4}
Among all various classes of best rational approximations, the diagonal sequences $r\in{\cal R}(k,k)$ of the {\em Walsh 
table} of $t^\alpha$, $\alpha\in(0,1)$ are studied in greatest detail \cite{Newman64,Ganelius79,Stahl93,SS93}. There 
is an existence and uniqueness of the BURA elements for all $k$ and $\alpha$. The distribution of poles, zeros, and 
extreme points of those elements plays a central role in asymptotic convergence analysis, when $k\to\infty$, 
thus is well known. In this section, we show that the above diagonal class perfectly fits within our positive 
$\calA^{-\alpha}$ approximation framework. 

{First, we collect some preliminary results that will be later needed for the proof of the main theorem. The 
following characterization lemma, which we state here without proof, is vital for our further investigations.
{\lemma\cite[Lemma 2.1]{SS93}\label{lemma1} {Let $m=k$ and $0<\alpha<1$. 
\begin{itemize}
\item[(a)] The best rational approximant $r^1_\alpha$ is of exact numerator and denominator degree $k$.
\item[(b)] All $k$ zeros $\zeta_1,\dots,\zeta_k$ and poles $d_1,\dots,d_k$ of $r^1_\alpha$ lie on the negative 
half-axis $\RR_{<0}$ and are interlacing; i.e., with an appropriate numbering we have
\begin{equation}\label{eq:zeros and poles}
0>\zeta_1>d_1>\zeta_2>d_2>\dots>\zeta_k>d_k>-\infty 
\end{equation}
\item[(c)] The error function $t^{1-\alpha}-r^1_\alpha(t)$ has exactly $2k+2$ extreme points $\eta_1,\dots,\eta_{2k+2}$ 
on $[0,1]$, and with an appropriate numbering we have
\begin{align}\label{eq:eta1}
0&=\eta_1<\eta_2<\dots<\eta_{2k+2}=1\\\label{eq:eta2}
\eta^{1-\alpha}_j-r^1_\alpha(\eta_j)&=(-1)^j E_\alpha(k,k;1),\qquad j=1,\dots,2k+2.
\end{align}
\end{itemize}
}}
The next lemma builds a bridge between the fractional decompositions of $r^1_\alpha$ and $t^{-1}r^1_\alpha$.
{\lemma\label{lemma2} {Let $m=k$, $0<\alpha<1$, and 
$$r^1_\alpha(t)=b^\ast_0+\sum_{j=1}^k\frac{c^\ast_j}{t-d_j},\qquad
t^{-1}r^1_\alpha(t)=\frac{c_{0,1}}{t}+\sum_{j=1}^k\frac{c_j}{t-d_j}.$$
Then
\begin{equation}\label{eq:coeffs correspondence}
c_{0,1}= E_\alpha(k,k;1),\qquad c_j=c^\ast_j/d_j,\quad j=1,\dots,k.
\end{equation}
}}
\begin{proof}
The second part of \eqref{eq:coeffs correspondence} follows directly from
$$\frac{1}{t(t-d_j)}=\frac{1}{d_j}\left(\frac{1}{t-d_j}-\frac{1}{t}\right),\quad j=1,\dots,k.$$
For the first part, we combine the above identity with \eqref{eq:eta1} and \eqref{eq:eta2}
$$c_{0,1}=b^\ast_0-\sum_{j=1}^k\frac{c^\ast_j}{d^j}=r^1_\alpha(0)=-\big(\eta_1^{1-\alpha}
-r^1_\alpha(\eta_1)\big)=E_\alpha(k,k;1).$$
The proof of the lemma is completed.
\end{proof}
Our last lemma provides an asymptotic bound on $E_\alpha(k,k;1)$. The proof can be found in \cite{Stahl93}.
{\lemma\cite[Theorem 1]{Stahl93}\label{lemma3} {The limit $$\lim_{k\to\infty}e^{2\pi\sqrt{\alpha 
k}}E_{1-\alpha}(k,k;1)=4^{1+\alpha}|\sin{\pi\alpha}|$$ holds true for each $\alpha>0$.}}

Now, we are ready to formulate and prove our main result.
{\theorem\label{thm1} {Let $\beta=1$ and $m=k$. For every normalized SPD M-matrix $\calA$ and every $\alpha\in(0,1)$, 
the matrix $\calA^{-1}r^1_\alpha(\calA)$ is doubly nonnegative and for all $\gamma\in\RR$
\begin{equation}\label{eq:thm1}
\frac{\|\calA^{-1}r^1_\alpha(\calA)\bff-\calA^{-\alpha}\bff\|_{\calA^{\gamma+1}}}{\|\bff\|_{\calA^{\gamma-1}}}
\le 4^{2-\alpha}|\sin{\pi(1-\alpha)}|e^{-2\pi\sqrt{(1-\alpha) k}}(1+\mathrm o(1)).
\end{equation}
}}
\begin{proof}

Based on the results in Lemma~\ref{lemma1}, we can quickly derive $\calA^{-1}r^1_\alpha(\calA)\ge0$. For this 
purpose, we study the sign pattern of $c$ and $d$ and assure the applicability of Proposition~\ref{prop3}. From 
\eqref{eq:zeros and poles} we know that all the poles $\{d_j\}$ are real, negative, and of multiplicity 1. The same 
holds true for the zeros $\{\zeta_j\}$. Since $r^1_\alpha(t)$ is continuous on $\RR\setminus\{d_j\}$, the function 
changes its sign $2k$ times - at each zero $\zeta_j$ and at each pole $d_j$. In Lemma~\ref{lemma2}, we have already 
computed that $r^1_\alpha(0)=c_{0,1}=E_\alpha(k,k;1)>0$, thus, due to interlacing, at each pole $d_j$ we have
$$\begin{array}{c} \lim_{t\to d_j^+}r^1_\alpha(t)<0\\ \lim_{t\to d_j^-}r^1_\alpha(t)>0\end{array}\quad\Longrightarrow
\quad\begin{array}{c} \lim_{t\to d_j^+}r^1_\alpha(t)=-\infty\\ \lim_{t\to d_j^-}r^1_\alpha(t)=+\infty\end{array}
\quad\Longrightarrow\quad c^\ast_j<0.$$
Since $c^\ast <0$ and $d<0$, from Lemma~\ref{lemma2} it follows that $c_0>0$ and $c>0$. Hence, Proposition~\ref{prop3} 
gives rise to $\calA^{-1}r^1_\alpha(\calA)\ge0$.

The error estimate \eqref{eq:thm1} is a direct corollary of Proposition~\ref{prop2} and Lemma~\ref{lemma3}. 
\end{proof}
}
Some remarks are in order. For any fixed $\alpha<1$, the relative error \eqref{eq:thm1} decays exponentially as
$k\to+\infty$ with order $\sqrt{(1-\alpha)k}$. When $\alpha\to 1$ the relative error decays linearly independently of
$k$, since $|\sin{\pi(1-\alpha)}|\to\pi(1-\alpha)$. It is straightforward to extend the coefficient
correspondence \eqref{eq:coeffs correspondence} to $c_j=c^\ast_j/d^\beta_j$ for any $(m,k,\beta)$,
such that $m<k+\beta$. Therefore, a necessary condition for Proposition~\ref{prop3} to be applicable is the sequence
$c^\ast$ to have constant sign. Due to the proof of Theorem~\ref{thm1}, it implies that zeros $\{\zeta_i\}_1^m$ and
poles $\{d_j\}_1^k$ of $r^\beta_\alpha$ should be interlacing, thus $|m-k|\le1$. Furthermore (see \cite[(20)]{HLMMV16})
the following identity always holds true
$$c_{0,1}+\sum_{j=1}^k c_j=0,\qquad m<k+\beta-1.$$
Hence, another necessary condition for applicability of Proposition~\ref{prop3} is $m\ge k+\beta-1$. Combining all
derived constraints, we observe that $\calA^{-\beta}r^\beta_\alpha(\calA)$ could be represented as a sum of doubly
nonnegative matrices only if $m=k+\beta-1$ and $m\le k+1$, meaning that we are left with the admissible triples
$$(m,k,\beta)= \{(k,k,1),(k,k,2),(k+1,k,2)\}.$$
In \cite{SS93} it is remarked that for the case $(k,k,2)$ it cannot be theoretically excluded that one root $\zeta_i$
and one pole $d_i$ of the 2-BURA $r^2_\alpha(t)$ lie outside of $\RR_{<0}$. For the case $(k+1,k,2)$ we can show that
$c_{0,2}=-E_\alpha(k+1,k;\beta)<0$, thus the assumptions of the proposition are again violated. In conclusion, the
triple $(k,k,1)$, investigated in Theorem~\ref{thm1} is the unique choice of parameters for which one can prove
positiveness of the approximation $\calA^{-\beta}r^\beta_\alpha(\calA)$ using Proposition~\ref{prop3}.

\section{Numerical tests}\label{section5}
{
The main goal of the numerical tests is to illustrate the positive properties
of the proposed approximations of inverse of fractional powers of 
SPD M-matrices. Complementary, we provide a short discussion related to some 
interpretations of the results from Section \ref{section4} in the case of
fractional diffusion problems.     

The test problems are in 1D. In this case we are able easier to
compute the exact solutions of linear systems with fractional powers 
of the corresponding tridiagonal matrix. 
Note that this does not cause restrictions to the derived conclussions.
As it was shown in \cite{HLMMV16}, if  $\Omega\subset\RR^d, d>1$, 
some PCG solver of optimal complexity (e.g., BoomerAMG)
can be utilized for efficient solution of the arising sparse linear systems, 
fully preserving the accuracy and efficiency of the composite 
algorithm.

We consider fractional powers of the Poisson's equation on the unit interval 
with Dirichlet boundary conditions:
\begin{equation}\label{eq:Poisson}
\mathcal L u:=-u''(x)=f(x),\qquad t\in[0,1],\quad u(0)=u(1)=0.
\end{equation}
}
On a uniform grid with mesh parameter $h=1/(N+1)$, using central 
finite differences, the operator $\mathcal L$ is 
approximated by the $N\times N$ matrix $\calA_h:=tridiag(-1,2,-1)/h^2$, 
which in turn can be rewritten as 
\begin{equation}\label{eq: normalized Laplacian}
\calA_h=4h^{-2}\calA, \qquad \calA:=tridiag\left(-\frac14, \frac12, -\frac14 \right). 
\end{equation}
The matrix $\calA$ is a normalized, SPD M-matrix, which eigenvectors and eigenvalues are explicitly known:
$$\Lambda_i=\sin^2\left(\frac{i\pi}{2(N+1)}\right),\qquad \bPsi_i=\left\{\sin\frac{im\pi}{N+1}\right\}_{m=1}^N,\qquad 
i=1,\dots,N.$$
We approximate $\mathcal L^\alpha$ by $\calA_h^\alpha$ and, due to Theorem~\ref{thm1}, the $\ell^2$ relative error is bounded by an $h$-dependent constant.
{\corollary\label{cor1} {Let $\bfu_h:=\calA_h^{-\alpha}\bff=\left(\frac{h}{2}\right)^{2\alpha}\calA^{-\alpha}\bff$ be 
the exact solution of the discretized fractional Poisson's equation $\calA_h^\alpha\bfu=\bff$. Let 
$r^1_\alpha\in\mathcal R(k,k)$ be 1-BURA and denote by 
$\bfu_{h,r}:=\left(\frac{h}{2}\right)^{2\alpha}\calA^{-1}r^1_\alpha(\calA)\bff$. Then 
$$\frac{\|\bfu_{h,r}-\bfu_h\|_2}{\|\bff\|_2}\le\left(\frac 4 
h\right)^{2(1-\alpha)}|\sin{\pi(1-\alpha)}|e^{-2\pi\sqrt{(1-\alpha) k}}(1+\mathrm o(1)).$$}}
\begin{proof}
Indeed, let $\bfu=\calA^{-\alpha}\bff$ and $\bfu_r=\calA^{-1}r^1_\alpha(\calA)\bff$. Applying 
$$\|\cdot\|_2=\|\cdot\|_{\calA^0}\le\mathrm k(\calA)\|\cdot\|_{\calA^2}<h^{-2}\|\cdot\|_{\calA^2},$$
where $\mathrm k(\calA)$ is the condition number of $\calA$, we derive
\begin{equation}\label{eq:scale}
\frac{\|\bfu_{h,r}-\bfu_h\|_2}{\|\bff\|_2}\le h^{-2}\frac{\|\bfu_{h,r}-\bfu_h\|_{\calA^2}}{\|\bff\|_{\calA^0}}
\le h^{-2}\left(\frac{h}{2}\right)^{2\alpha}\frac{\|\bfu_r-\bfu\|_{\calA^2}}{\|\bff\|_{\calA^0}}.
\end{equation}
The result follows from \eqref{eq:thm1} for $\gamma=1$.
\end{proof}
Due to Corollary~\ref{cor1}, we can compute the minimal degree $k$ that guarantees $\|\bfu_{h,r}-\bfu_h\|_2/\|\bff\|_2 
<\varepsilon$ for every given pair $(\varepsilon,h)$. Such an $\ell^2$ error analysis is outside of the scope of this 
paper, so we will not further elaborate on it.   

\begin{table}
\caption{ Errors $E_\alpha(k,k;1)$ of BURA $r^1_\alpha(t)$ of $ t^{1 -\alpha}$ on $[0,1]$.}\label{t:error}
\centering
\begin{tabular}{|c|c|c|c|}
\hline
$\alpha$&$E_\alpha(5,5,1)$&$E_\alpha(6,6,1)$&$E_\alpha(7,7,1)$\\ \hline
  0.25 & 2.8676e-5& 9.2522e-6& 3.2566e-6\\
  0.50 & 2.6896e-4& 1.0747e-4& 4.6037e-5\\
  0.75 & 2.7162e-3& 1.4312e-3& 7.8966e-4\\
\hline
\end{tabular}
\end{table}
In our numerical experiments, we choose $\alpha\in\{0.25,0.5,0.75\}$ and $k\in\{5,6,7\}$. The square root of an M-matrix is 
again an M-matrix \cite{Alefeld82} and it is easy to check that $\calA_h^3$ is also an M-matrix. Therefore, all 
considered $\calA_h^{\alpha}$ are M-matrices, their inverse matrices are doubly nonnegative (but dense!), and 
constructing computationally cheep approximants within the same class is of great practical importance. The 
univariate error estimates $E_\alpha(k,k;1)$ for the above choice of parameters are summarized in Table~\ref{t:error}. 
Note that each of them satisfies the inequality \eqref{eq:thm1} even without introducing the low-order term $\mathrm 
o(1)$ in the right-hand-side. 
{
A modified Remez algorithm is used for the derivation of $r^1_\alpha$ 
\cite{PGMASA1987,CheneyPowell1987}.  
}

\begin{figure}
\centering
\includegraphics[width=0.32\textwidth]{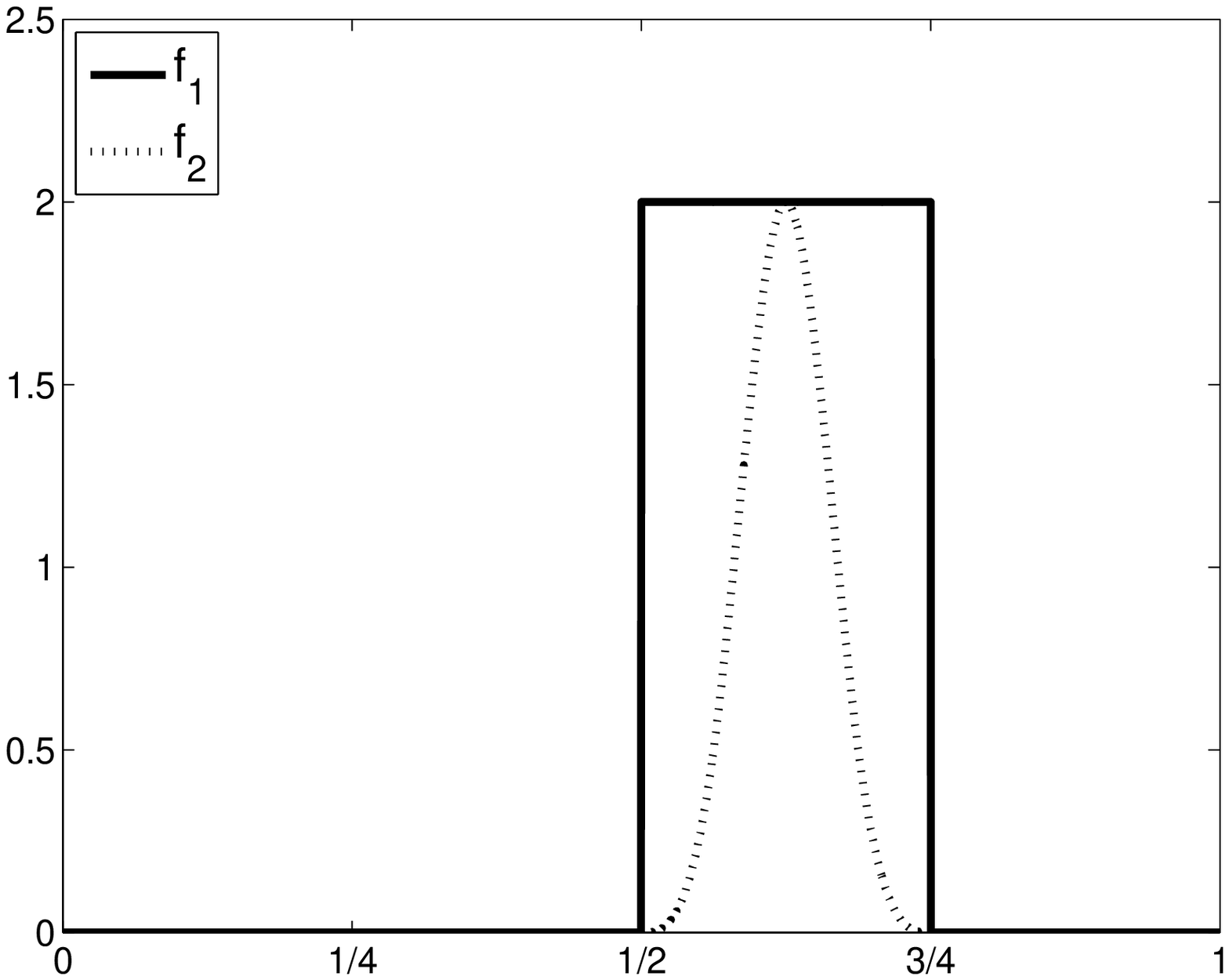}
\includegraphics[width=0.32\textwidth]{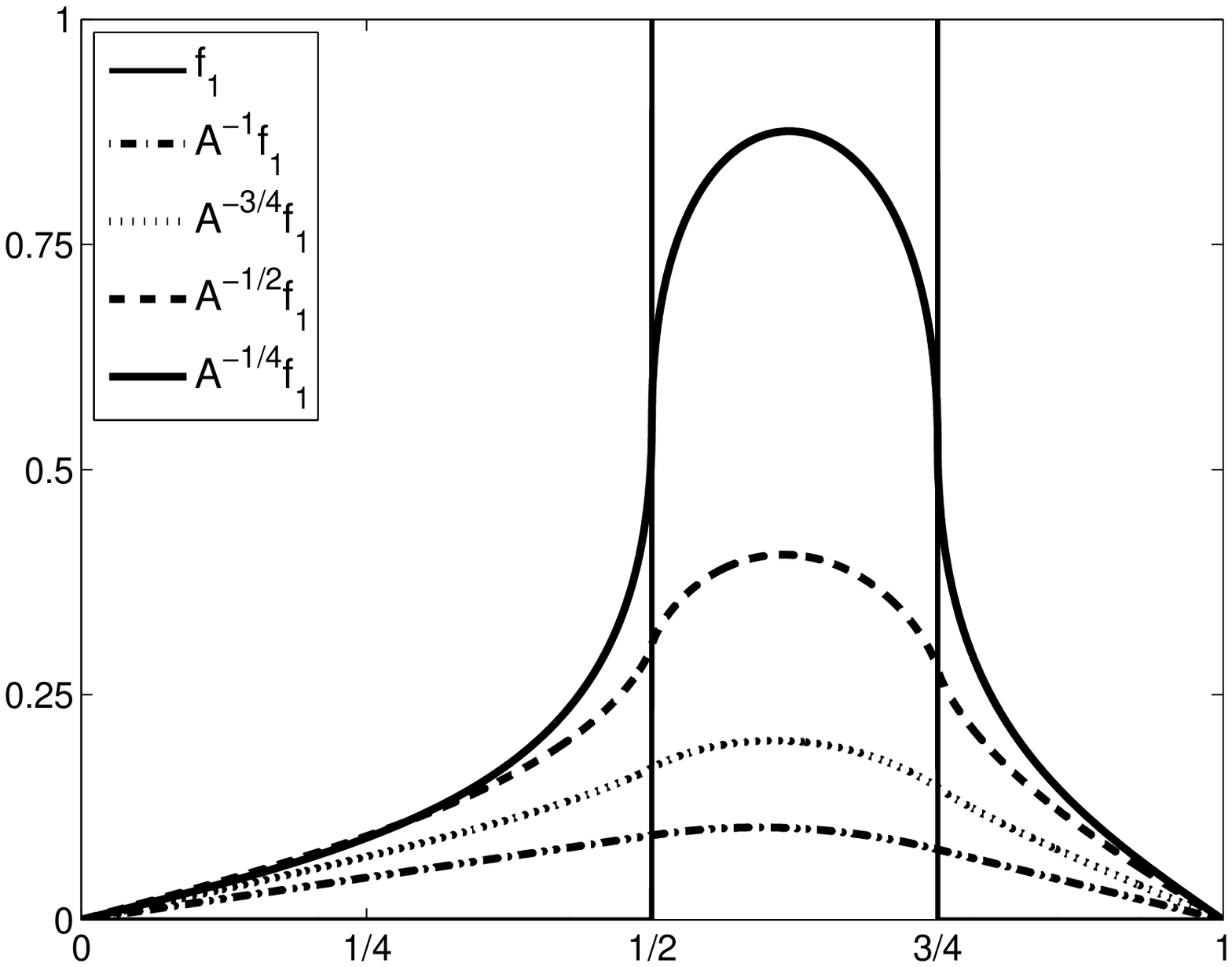}
\includegraphics[width=0.32\textwidth]{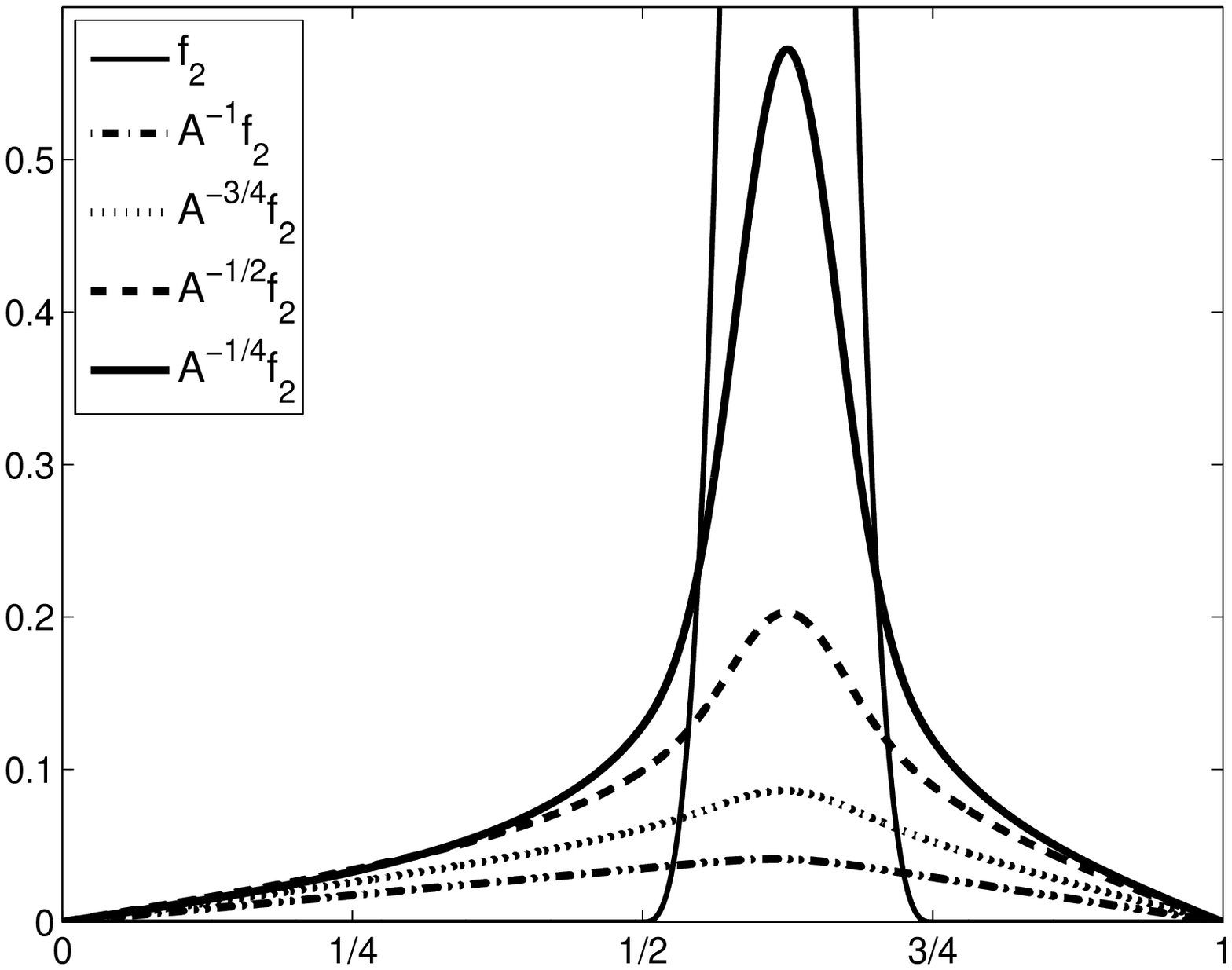}
\caption{Test data and their exact fractional diffusions.}\label{fig1}
\end{figure}

For $f$ in \eqref{NonLP} we take two different 
{positive}
 functions, supported on the interval $[1/2,3/4]$. The first one $f_1$ is 
piecewise constant and discontinuous, while the second one $f_2$ is a $C^2$ 
cubic spline function, corresponding to the Irwin-Hall distribution. Together with 
the exact discretized solutions $\calA_h^{-\alpha}\bff$, they are illustrated on 
Fig.~\ref{fig1}.

\begin{figure}
\centering
\includegraphics[width=0.32\textwidth]{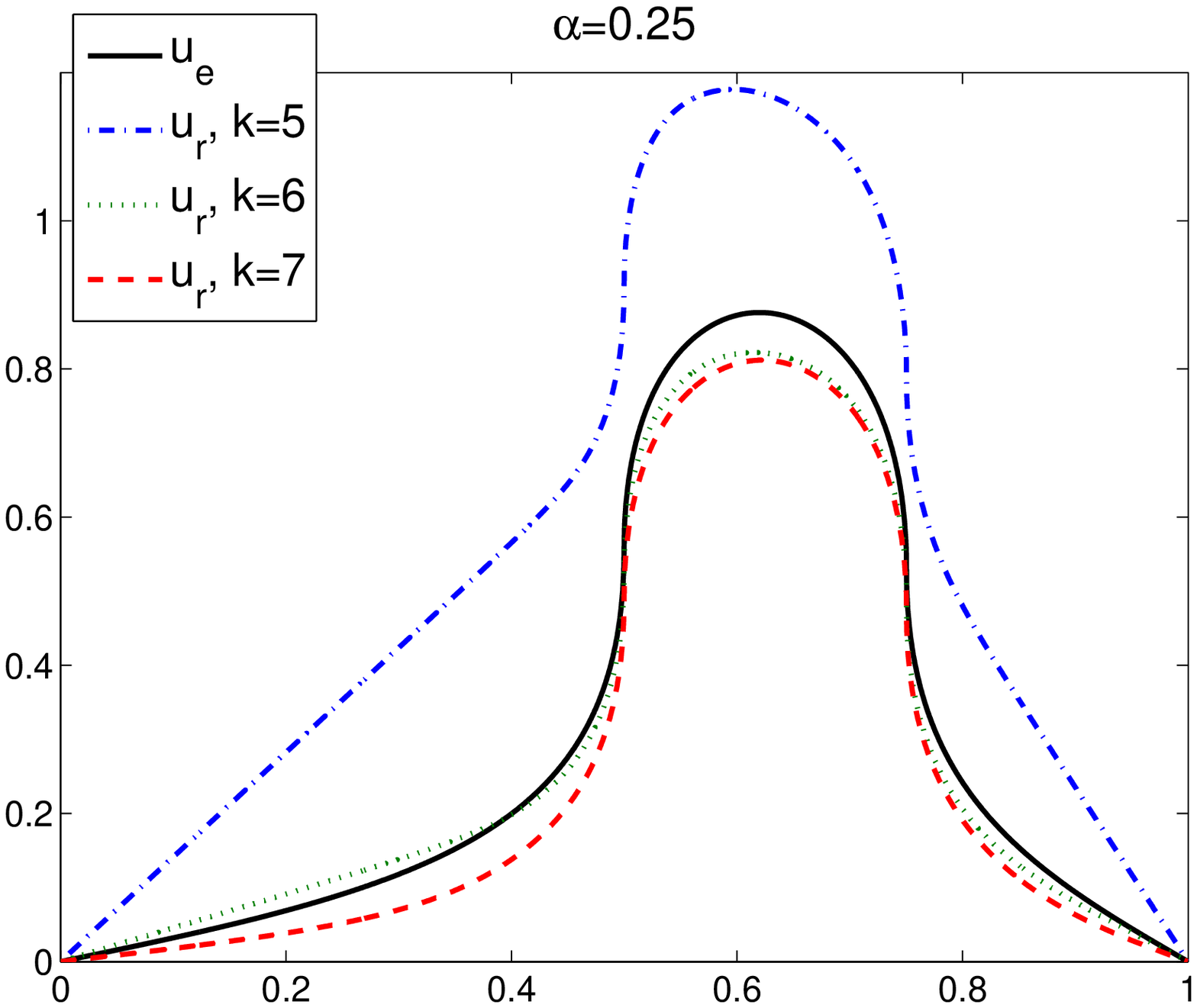}
\includegraphics[width=0.32\textwidth]{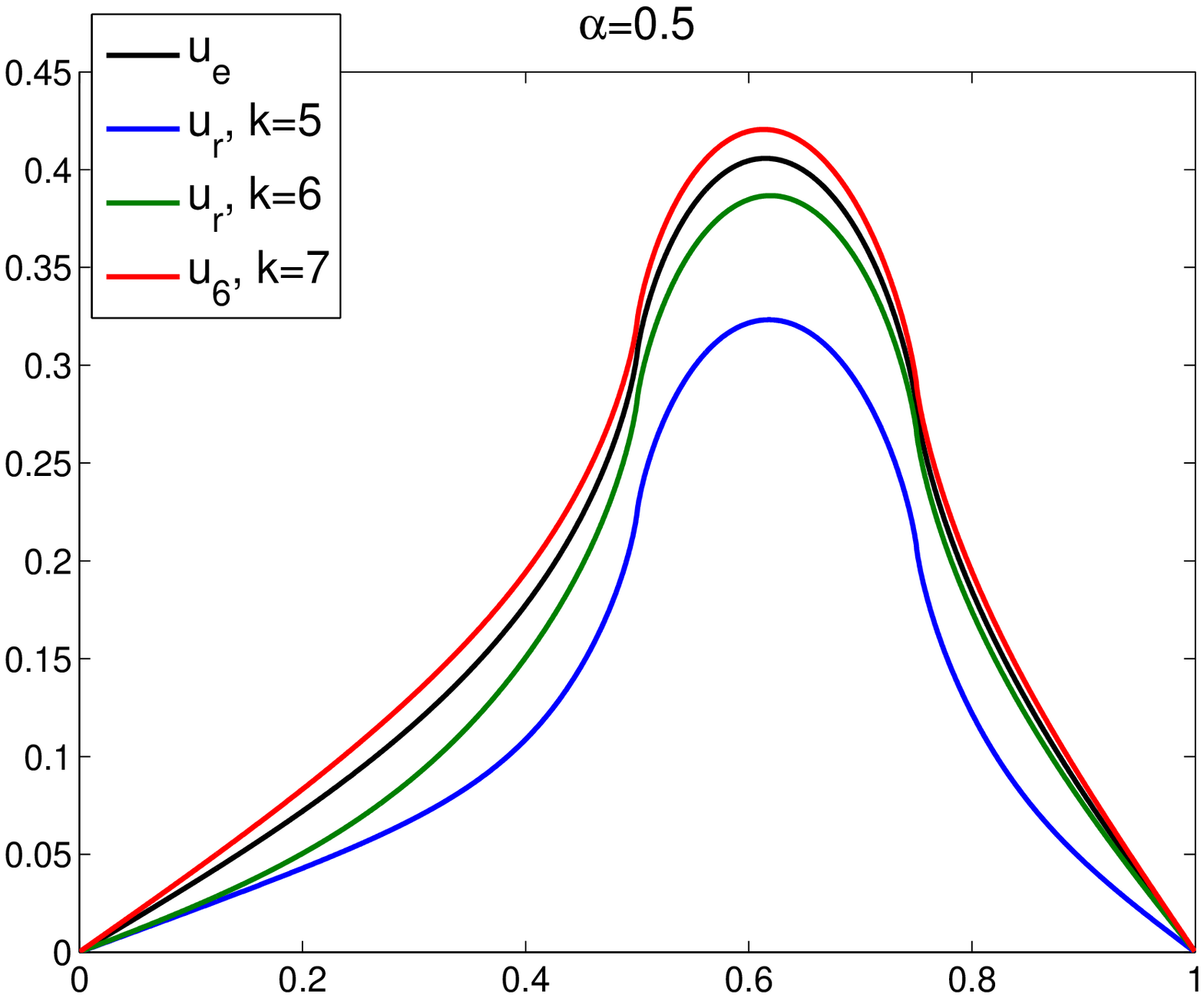}
\includegraphics[width=0.32\textwidth]{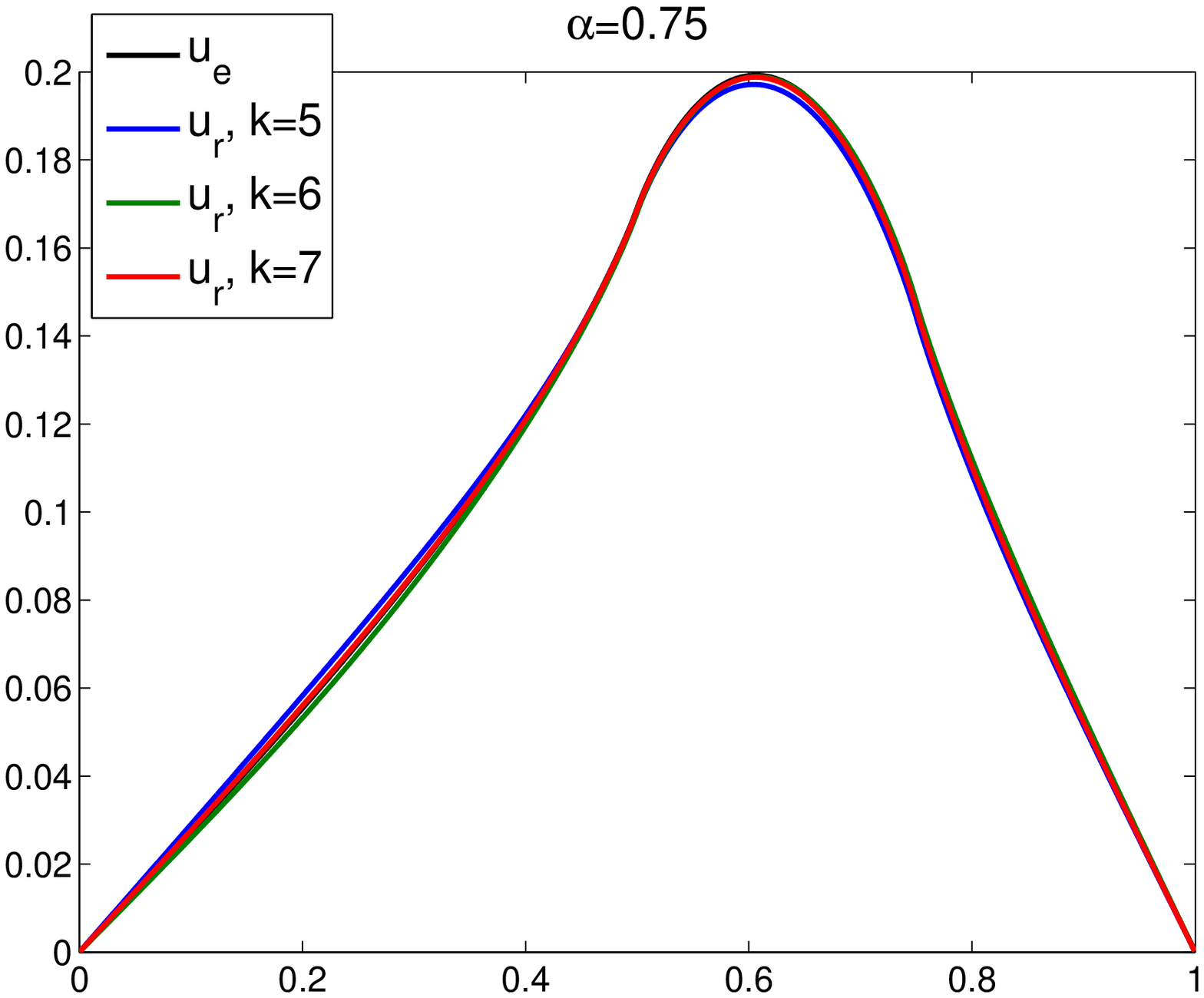}
\caption{Positive approximations of $\calA_h^{-\alpha}$ on $\bff_1$ for $h=2^{-11}$.}\label{fig2}
\end{figure}
\begin{figure}
\centering
\includegraphics[width=0.32\textwidth]{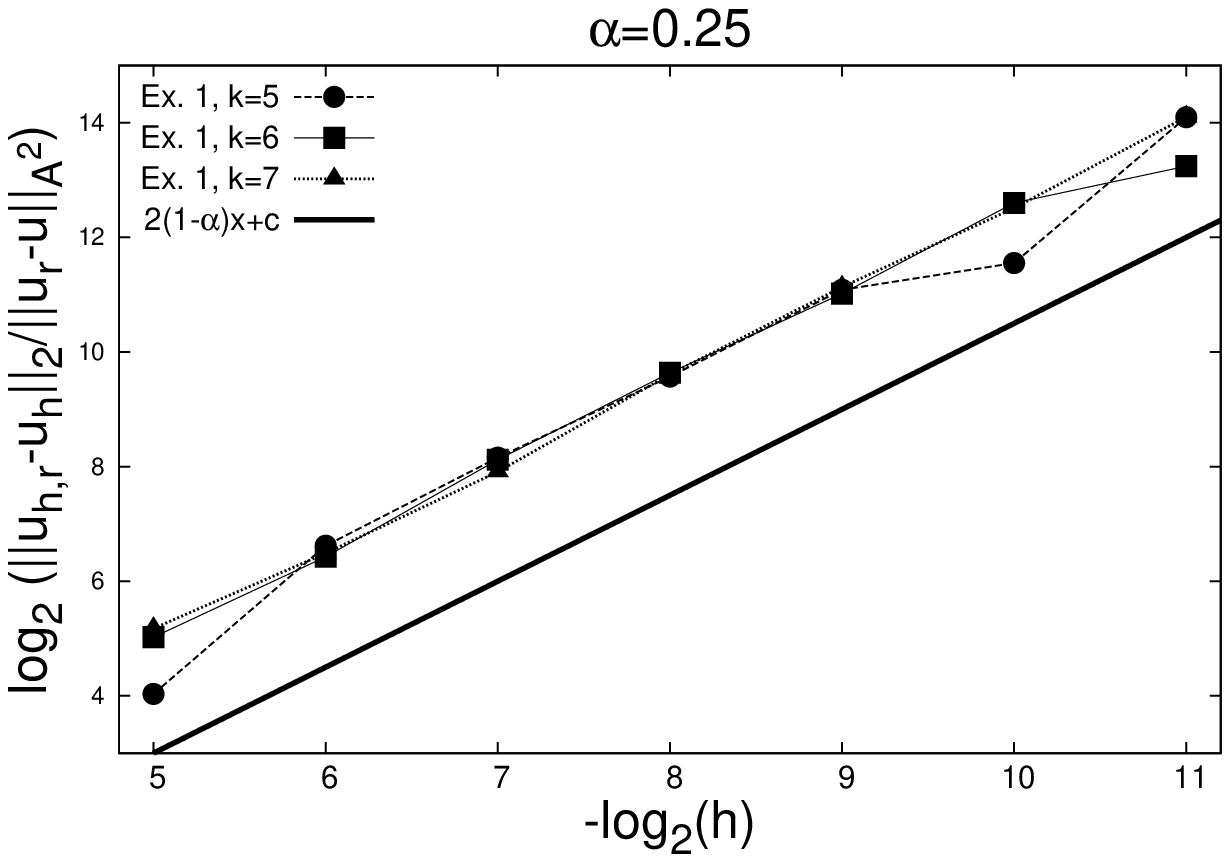}
\includegraphics[width=0.32\textwidth]{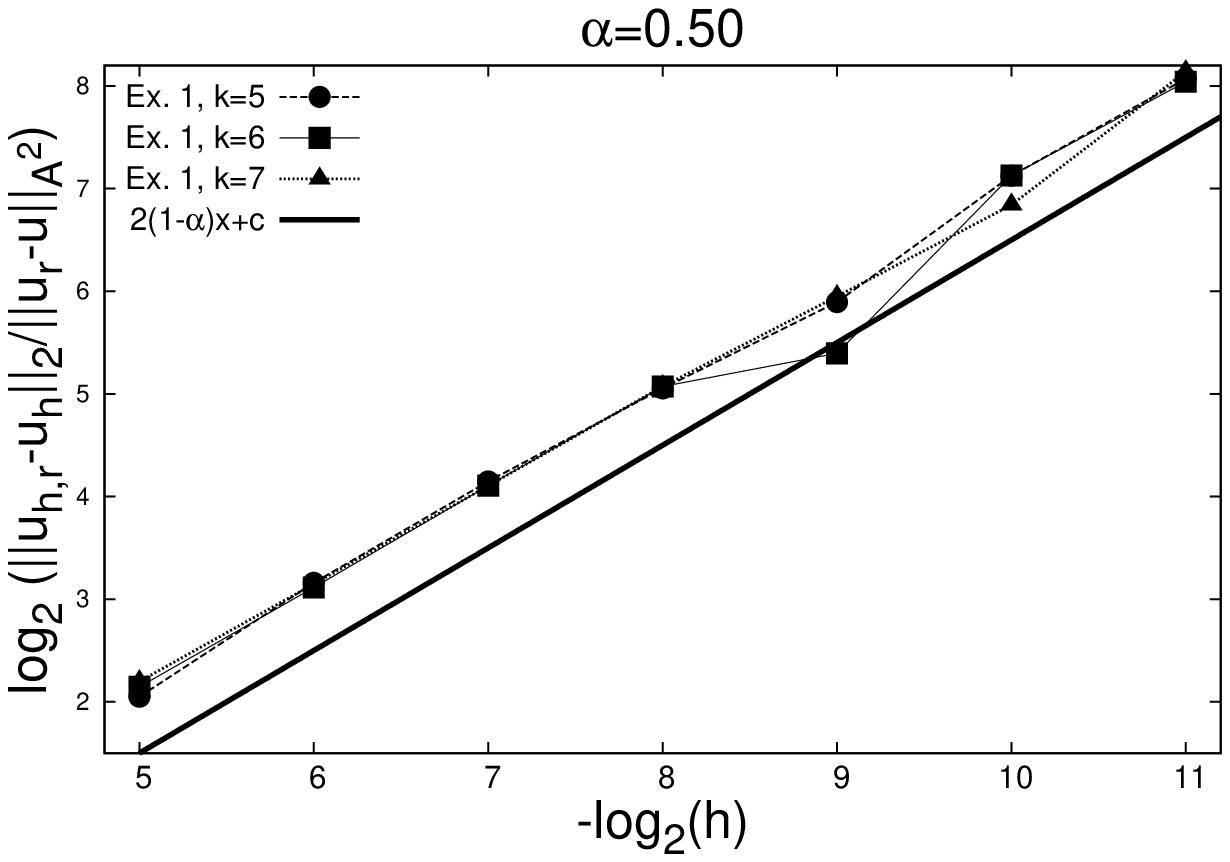}
\includegraphics[width=0.32\textwidth]{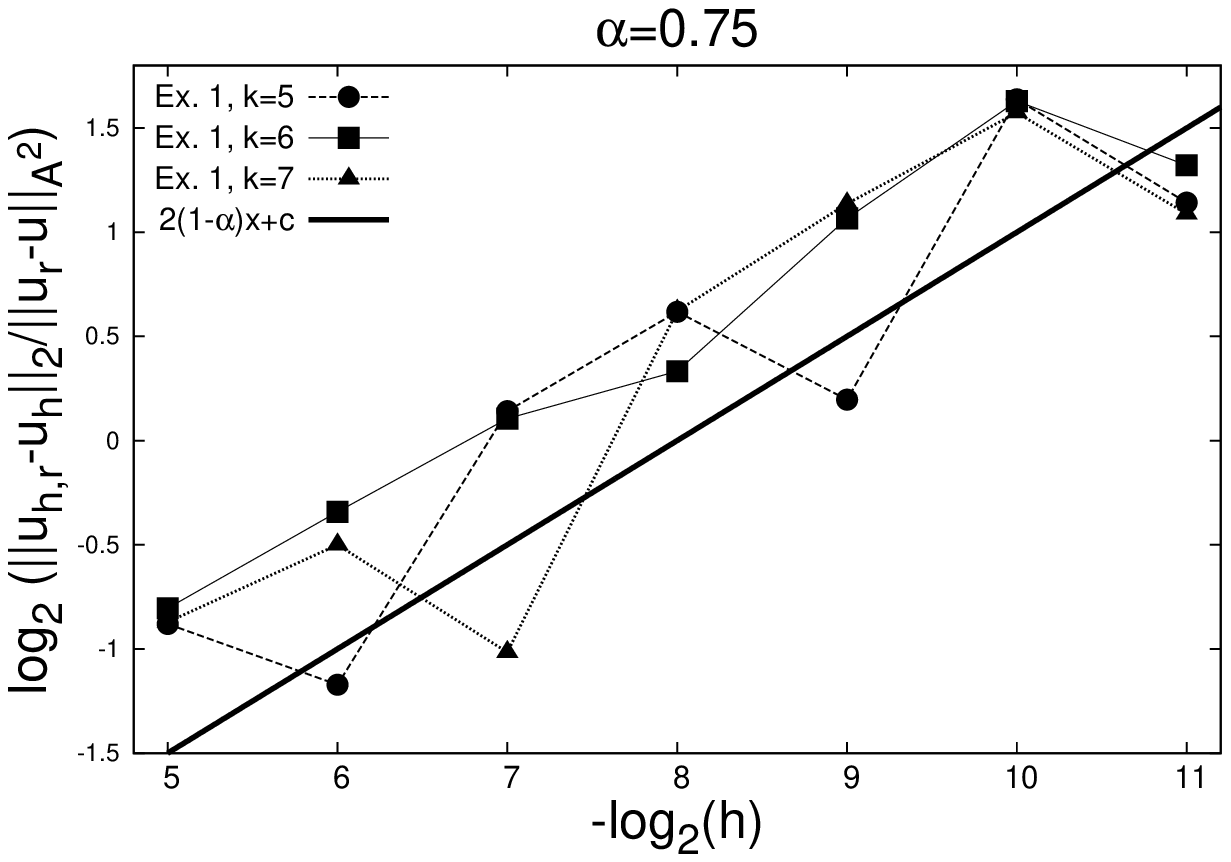}
\caption{Numerical confirmation for 
$\frac{\|\bfu_{h,r}-\bfu_h\|_2}{\|\bff_1\|_2}/\frac{\|\bfu_r-\bfu\|_{\calA^2}}{\|\bff_1\|_{\calA^0}}=\mathrm 
O(h^{-2(1-\alpha)})$ in \eqref{eq:scale}.}\label{fig3}
\end{figure}

We consider 
{mesh parameters}
 $h=2^{-m}$, $m\in\{5,6,\dots,11\}$. On Fig.~\ref{fig2} the corresponding approximants 
$\bfu_{h,r}$, $h=2^{-11}$, of $\bfu_h=\calA_h^{-\alpha}\bff_1$ are plotted. As suggested by Corollary~\ref{cor1}, 
$\bfu_{h,r}$ fails 
{
to approximate well 
$\bfu_h$ on a fine grid for smaller $k$ 
}
and $\alpha$ (see $k=5$, $\alpha\in\{0.25,0.5\}$). For 
{
larger $\alpha$, 
}
the exponential 
growth of the $\ell^2$ relative error with $h\to 0$ is less 
significant, as it is of order $2(1-\alpha)$, thus when $\alpha=0.75$ all the three approximants, corresponding to 
$k=\{5,6,7\}$ follow closely the graph of $\calA_h^{-3/4}\bff_1$. 
On Fig.~\ref{fig3} and in Table~\ref{table2} we numerically confirm the asymptotic behavior of the relative $\ell^2$ 
error from Corollary~\ref{cor1}. 

\begin{table}
\centering
\caption{{$\ell^2$ relative error $\frac{\|\bfu_{h,r}-\bfu_h\|_2}{\|\bff_2\|_2}$.}}\label{table2}
\begin{tabular}{|c|ccc|ccc|ccc|}
\hline
\multirow{2}{*}{$h$} & \multicolumn{3}{|c|}{$\alpha=0.25$} & \multicolumn{3}{|c|}{$\alpha=0.5$} & 
\multicolumn{3}{|c|}{$\alpha=0.75$}\\ \cline{2-10}
& $k=5$\hspace{1ex} & $k=6$\hspace{1ex} & $k=7$\hspace{1ex} & $k=5$\hspace{1ex} & $k=6$\hspace{1ex} & $k=7$\hspace{1ex} 
& $k=5$\hspace{1ex} & $k=6$\hspace{1ex} & $k=7$\hspace{1ex}\\ \hline
$2^{-5}$  & 6.3e-5 & 1.2e-4 & 7.5e-5 & 1.9e-4 & 3.1e-4 & 1.0e-4 & 8.5e-4 & 3.5e-4 & 2.1e-4 \\
$2^{-6}$  & 6.8e-4 & 2.1e-4 & 1.2e-4 & 9.3e-4 & 6.2e-4 & 2.6e-4 & 3.0e-4 & 7.3e-4 & 1.9e-4 \\
$2^{-7}$  & 4.9e-3 & 9.9e-4 & 2.2e-4 & 3.2e-3 & 1.3e-3 & 5.5e-4 & 1.6e-3 & 1.0e-3 & 6.6e-5 \\
$2^{-8}$  & 1.2e-2 & 4.9e-3 & 1.0e-3 & 5.0e-3 & 2.4e-3 & 1.0e-3 & 2.7e-3 & 6.7e-4 & 4.4e-4 \\
$2^{-9}$  & 3.7e-2 & 9.6e-3 & 4.9e-3 & 5.6e-3 & 1.3e-3 & 1.2e-3 & 8.2e-4 & 1.3e-3 & 1.1e-3 \\
$2^{-10}$ & 2.1e-2 & 3.6e-2 & 9.5e-3 & 2.4e-2 & 8.8e-3 & 1.4e-3 & 5.8e-3 & 3.0e-3 & 1.5e-3 \\
$2^{-11}$ & 1.9e-1 & 2.3e-2 & 3.6e-2 & 4.2e-2 & 1.4e-2 & 8.9e-3 & 1.6e-3 & 9.9e-4 & 4.3e-4 \\\hline
\end{tabular}

\end{table}

{
The sufficient conditions from Proposition \ref{prop3} hold true for the cases under
consideration. This means that positivity of all considered approximations is guaranteed.
Therefore, the discrete maximum principle is always inherited. The presented numerical 
results are fully aligned with the theory. What is very important is the numerical
robustness of positivity with respect to both accuracy parameters $h$ and $k$
which  is confirmed for all $\alpha \in \{0.25,0.5,0.75\}$. Even in the case of lower 
accuracy, we do not observe any oscillations. The monotonicity preservation of the data
is clearly expressed, capturing their geometrical shape.
}

\section{Concluding remarks}\label{section6}
This study is inspired by some quite recent results in the numerical 
methods for fractional diffusion problems. In the Introduction, we 
discussed three methods based on reformulation of the original nonlocal 
problem into local (elliptic, pseudo parabolic, and integral) problems. 
In all cases, the cost is in the increased dimension of computational domain 
from $d$ to $d+1$.  

Our approach is based on best uniform rational approximations of 
$t^{\beta - \alpha}$, $0\le t \le 1$. The primal motivation is to reduce 
the computational complexity. A next important step is made in this paper. 
Here, we provide sufficient conditions to guarantee positive approximation 
of the inverse of fractional powers of normalized SPD M-matrices. Therefore, we 
get a numerical method which preserves the maximum principle. The presented 
numerical results clearly confirm the monotone behaviour of the 
solution, without any observed oscillations. Further research has to be 
devoted to the topic of accuracy of mass conservation.

The currently available methods and algorithms for numerical solution of 
boundary value problems with  fractional power of elliptic operators 
have a quite different nature. A serious theoretical and experimental 
study is required to get a comparative analysis of their advantages and 
disadvantages for particular classes of problems. For instance, the error 
analysis is in different functional spaces assuming different conditions 
for smoothness. The comparison of the computational complexity is also an 
open question. 

As a part of our analysis, Theorem 1 provides a sharp error estimates 
for the 1-BURA based approximations $E(k,k,1)$. Then, at the beginning 
of Section 5, we showed how this result can be used to derive relative 
error estimates of the numerical solution of fractional order elliptic 
problems in $\ell^2$. The numerical tests are well aligned with this 
theoretical estimates. The presented approach has a strong potential 
for further development addressing different pairs of functional spaces in 
the relative error estimates, varying the smoothness assumptions, 
for $d=1,2,3$. 

In addition, a lot of new numerical tests are needed to 
evaluate/confirm/compare the computational efficiency for more 
realistic towards real-life large-scale super diffusion problems.
\bibliographystyle{abbrv}

\end{document}